\documentclass[12pt]{article}

\usepackage{amsfonts}
\usepackage{amsmath}
\usepackage{graphicx}

\newcommand{\pibar}{\overline{\Pi}}

\newcommand{\Xbar}{\overline{X}}

\newcommand{\topr}{\stackrel{\mathrm{P}}{\longrightarrow}}

\newcommand{\todr}{\stackrel{\mathrm{D}}{\longrightarrow}}

\newcommand{\R}{\Bbb{R}}

\newcommand{\rmd}{{\rm d}}

\newcommand{\dto}{\downarrow}

\newcommand{\be}{\begin{equation}}
\newcommand{\ee}{\end{equation}}
\newcommand{\bea}{\begin{eqnarray}}
\newcommand{\eea}{\end{eqnarray}}
\newcommand{\bean}{\begin{eqnarray*}}
\newcommand{\eean}{\end{eqnarray*}}
\newcommand{\ben}{\begin{equation*}}
\newcommand{\een}{\end{equation*}}
\newcommand{\ba}{\begin{aligned}}
\newcommand{\ea}{\end{aligned}}

\def\nexto{\kern -0.54em}

\def\topr{\buildrel P \over \to }

\numberwithin{equation}{section}
\numberwithin{theorem}{section}
\numberwithin{corollary}{section}
\numberwithin{proposition}{section}
\numberwithin{lemma}{section}

\usepackage{authblk}
\author[1]{Ross A. Maller\thanks{ross.maller@anu.edu.au}}
\author[2]{Yuguang Fan\thanks{yuguang.fan@unimelb.edu.au}}
\affil[1]{ \small Research School of Finance,  Actuarial Studies and Applied Statistics, the Australian National University}
\affil[2]{ \small School of Mathematics and Statistics, the University of Melbourne}
\begin{document}

\title{\bf Thin and Thick Strip Passage Times for L\'evy Flights and L\'evy Processes}

\maketitle
\begin{abstract}
\noindent
We review some of the theory relevant to passage times of one-dimensional L\'evy processes out of bounded regions, highlighting results that are useful in physical phenomena modelled by heavy-tailed L\'evy flights. The process is hypothesised to describe the motion of a particle on the line, starting at 0, and exiting either a fixed interval $[-r,r]$, $r>0$,  or a time-dependent, expanding, set of intervals of the form  $[-rt^\kappa, rt^\kappa]$, $r>0$, $\kappa>0$.  
Asymptotic behaviour of the exit time  may be as $r\dto 0$ or as $r\to\infty$, but particular emphasis is placed herein on ``small time" approximations, corresponding to exits from or transmissions through thin strips.
Applications occur for example in the transmission of photons through moderately doped thin or thick wafers by means of ``photon recycling", and
in atmospheric radiation modelling.
\end{abstract}

\section{Introduction}\label{s0}
L\'evy flights are increasingly being proposed as models for physical phenomena in a variety of applications areas. Convincing empirical evidence for their existence is reported for example in experiments or observations on the transmission of photons through scattering media such as clouds, especially relevant to climatological studies (e.g., \cite{scatter}) and the transport of ``holes" in semiconductors (e.g., \cite{holes}).
The latter two applications, in particular, are concerned with the one-dimensional motion of  a particle or object in a bounded domain -- a ``slab" in the case of \cite{scatter}, a ``wafer" in the case of \cite{holes}, by means of jumps whose magnitude has a heavy tailed distribution, specifically, a power law tail. The corresponding random walk -- the cumulative sum of the jump sizes -- has an infinite variance
but is in the domain of attraction of a stable law.  Since stable (L\'evy) processes  in continuous time can be obtained as limits of interpolated L\'evy flights, we can expect the particles' motions to be well described by a stable process in one (or more) dimensions.
In many contexts these provide better descriptions of the process than a Brownian motion model.

In the present paper we review some of the theory relevant to passage times of one-dimensional L\'evy processes out of bounded regions with a view to highlighting results  useful in these kind of modelling contexts.
The results will be formulated in terms of general random walks and L\'evy processes taking values in the real line $\R$, but special case reductions to (one-dimensional) stable processes and L\'evy flights will be made along the way.

A process $X_t$, $t\ge 0$, starting at 0, is hypothesised to describe the motion of a particle on the line, starting at 0, and exiting either a fixed interval $[-r,r]$ or a time-dependent, expanding, set of intervals of the form
$[-rt^\kappa, rt^\kappa]$.   Here $\kappa\ge 0$ is a parameter to be specified and $r>0$ is a variable which in the physical context (and in the case $\kappa=0$) would specify the thickness of the slab or wafer. 
Introducing the time dependence allows for the slab thickness to increase in dimension at an algebraic rate in time, starting from zero thickness. Of course in the physical context this would represent a limiting case where the particle's escape is restricted in some way.

Exits of $X_t$ from the regions just described are represented by passage times of the form
\be\label{defT}
T_\kappa(r)=\inf\{t\ge 0:|X_t|>rt^\kappa\},\ r>0.
\ee
Expressed in space-time coordinates, 
 $T_\kappa(r)$ is the time at which the point $(t,X_t)$ first  falls outside the rectangular (when $\kappa=0$) or curvilinear (when $\kappa>0$) region defined by
 \be\label{reg}
 R_\kappa:=\{(t,x):t>0, x\in\R, |x|\le rt^\kappa\},
 \ee
having started at $(0,0)$.  When $\kappa=0$ this corresponds to exit from a finite slab of width $2r$.


There are two possible asymptotics of interest.
Most common in applications is to take $r\to\infty$, a ``large region", hence ``large time", approximation.
In \cite{scatter}, for example, the application is to photon free path distribution in atmospheric radiation modelling, where 
the heavy tailed distribution of path lengths arises with the passage of light through cloud types where rapid aerosol density fluctuations result in free paths varying from less than 1 meter to many kilometres. ``(L\'evy-stable) distributions are therefore well suited to model photon free paths globally" (\cite{scatter}).

At the other end of the scale, in \cite{holes}, short
wavelength radiation is shone onto doped crystal wafers of thicknesses of the order of 50-350 $\mu$m,  creating an electron-hole pair.  On a nanosecond time scale, the hole moves randomly with the thermal velocity
(this motion well described by Brownian motion)  until recombining with an electron. Emitted photons are absorbed, thereby generating a new hole, from which a new photon is emitted.
 The resultant ``photon recycling" propagates through the crystal and its location, when projected on the $y$-axis, defines  a 1-dimensional random walk. But in this instance the jump distribution is better described by a power law tail decay having infinite variance. After a large number of jumps a stable distribution is to be expected for the particle's location, and this is indeed observed experimentally.  ``The power law decay of the hole concentration $\ldots$ is steep enough at short distances $\ldots$ to fit the data for the thin sample,  and at the same time slow enough at large distances $\ldots$ to account for the data for thick samples" (\cite{holes}, p.20). 
Here, for thin samples,  concern is with the short time behaviour of the approximating stable process and the asymptotic of interest for  $T_\kappa(r)$ is as $r\dto 0$.

In the next sections we consider both ``thin strip" (short time) and large time behaviour.
First we need some preparatory material in the next section.

\section{Functional Laws for Random Walks}\label{s2} 
``Short Time" only makes sense for continuous time processes.
How to formulate these in a way that's relevant to the physical situation?
We can start with a  random walk $S_n$ comprised of i.i.d jumps $\xi_i$, $i=1,2,\ldots$; thus,
\be\label{defS}
S_n=\sum_{i=1}^n \xi_i,\ n=1,2, \ldots,
\ee
with $S_0=0$. Each $\xi_i$ has the distribution of a generic rv $\xi$, which is assumed to have a heavy tail satisfying
\be\label{tail}
P(|\xi|>x) \sim x^{-\alpha}, \ {\rm as}\ x\to\infty,
\ee
 as well as the tail ``balance" condition
\be\label{bal}
\lim_{x\to\infty}\frac{P(\xi>x)}{P(|\xi|>x)}\
{\rm exists\ and\ equals}\ c,\  0\le c\le 1.
\ee
In \eqref{tail},  $\alpha$ is a parameter in $(0,2)$.
Conditions \eqref{tail} and \eqref{bal} imply that the random walk is in the large-time domain of attraction of a stable random variable, $S^{(\alpha)}$, of index $\alpha$, which means that there are centering and norming  sequences  $A_n$ and $B_n$, such that 
\be\label{doa}
\frac{S_n-A_n}{B_n}\todr S^{(\alpha)},\ {\rm as}\ n\to\infty.
\ee
Here and below ``$\todr$" denotes convergence  in distribution.
When $0<\alpha<1$ in \eqref{tail}, the $A_n$ can be taken as 0, and when $1<\alpha<2$, $E|\xi|$ is finite and we can set $A_n=nE(\xi)$, which we can further take as 0 by recentering. Further, when \eqref{tail} and \eqref{bal} hold, the $B_n$ satisfy $B_n\sim an^{1/\alpha}$ as $n\to\infty$ for some constant $a>0$ which can be taken as 1 by rescaling.\footnote{When $\alpha=1$ the centering is logarithmic. More general norming sequences, where $n^{1/\alpha}$ is multiplied by a slowly varying function, rather than a constant, can be catered for in a more general setup. But the situations we consider here suffice to illustrate the main points.}

Assume then that we have the convergence
\be\label{dona}
\frac{S_n}{n^{1/\alpha}}\todr S^{(\alpha)},\ {\rm as}\ n\to\infty.
\ee
This extends to a ``functional" version as follows. For each $t$ in a finite interval $[0,T]$ let $S_n(t)$ be the interpolated process defined by ``joining the dots" in the graph $(j,S_j)$, $j=0,1,2,\ldots,n$. Then we have, in an appropriate sense, the convergence
\be\label{dofna}
\frac{S_n(t)}{n^{1/\alpha}}\todr S_t^{(\alpha)},\
{\rm as}\ n\to\infty,
\ee
where $(S_t^{(\alpha)})$, $t>0$, is a stable process of index $\alpha$.    



We can define discrete-time passage times by  
\be\label{defTS}
T^S_\kappa(r)=\min\{n=1,2,\ldots:|S_n|>rn^\kappa\},\ r>0,
\ee
 and in view of the functional law in \eqref{dofna}, they will  converge as $n\to\infty$, after rescaling,  to the corresponding passage times for $S_t^{(\alpha)}$, defined by
 \be\label{Talpha}
T^{(\alpha)}_\kappa(r)=\inf\{t>0:|S_t^{(\alpha)}|>rt^\kappa\},\ r>0.
\ee
Having adopted a L\'evy flight $(S_n)_{n=1,2,\ldots}$ or a stable process $(S_t^{(\alpha)})_{t>0}$ as a physical model, properties of the time taken for the process, starting at 0 at time 0, to traverse and exit the regions $R_\kappa$ defined in \eqref{reg}, can then be derived.
In the continuous time process we can analyse small time behaviour, as $t\dto 0$, and this translates for the passage time to small ``slab width" behaviour, i.e., as $r\dto 0$. Large time behaviour, as $r\to\infty$,  can similarly be considered.  For the discrete time model, only large $n$  behaviour can be considered.

There are two ways to think about the modelling procedure. One, is that we model the physical process by the discrete time random walk $S_n$ and think of the stable process as a continuous time approximation to it. Alternatively, we could simply model the physical process directly with the continuous time stable process. Much of the time the difference will just be one of interpretation because in practice the continuous time model has to be fitted via a discretised version in any case.

\section{Functional Laws for L\'evy Processes}\label{s3} 
Rather than the discrete time L\'evy flight, we may adopt as a physical model a general L\'evy process $(X_t)_{t\ge 0}$ 
and concern ourselves directly with distributions and behaviour as $r\dto 0$ or as $r\to\infty$ of the passage times $T_\kappa(r)$ defined in \eqref{defT}.
Just as for the discrete time random walk, a real-valued L\'evy process $(X_t)_{t\ge 0}$ may be in the domain of attraction of a stable law either as $t\dto 0$ or as $t\to\infty$. Conditions on the L\'evy tail characterising this are available in \cite{mm} 
 (the canonical triplets of the rescaled underlying processes are required to converge to that of the stable process in an obvious way). 
If these conditions obtain, there is a functional law  analogous to \eqref{dofna} which gives that a rescaled version of $X_t$ converges  to a stable process  $(S_t^{(\alpha)})_{t\ge 0}$
either as $t\dto 0$ or as $t\to\infty$.

\section{Basic Stable Process Properties}\label{s4} 
Stable processes, $S_t^{(\alpha)}$, $t\ge 0$, are pure jump L\'evy processes whose distributions at each fixed time $t>0$ are stable with index $\alpha$ in $(0,2)$. The  distributions are absolutely continuous; thus, the distributions of  $S_t^{(\alpha)}$ at each fixed time $t>0$ possess densities. Moments of $S_t^{(\alpha)}$ of orders less than $\alpha$ are finite but  moments of orders greater than or equal to $\alpha$ are infinite.   Particularly attractive for physical modelling is that stable processes have the {\it scaling property:}
for each $\lambda>0$, $(S_{t\lambda}^{(\alpha)})_{t\ge 0}$ has the same distribution, as a process, as 
$(\lambda^{1/\alpha} S_t^{(\alpha)})_{t\ge 0}$.
In particular, for each $\lambda>0$, $S_{\lambda}^{(\alpha)}$
has the same distribution as  $\lambda^{1/\alpha} S_1^{(\alpha)}$, so in a  sense which can be made precise, $S_{\lambda}^{(\alpha)}$ is of order  $\lambda^{1/\alpha}$ both for small and for large $\lambda$. Note also that the scale invariance implies a form of self-similarity of the process.

A L\'evy process is completely determined by its {\it canonical triplet} 
$(\gamma, \sigma^2, \Pi)$, where $\gamma\in\R$ is a shift constant, $\sigma^2\ge 0$ allows for a Brownian component (for the stable laws, $\sigma^2=0$), and $\Pi$ is a measure on $\R$ which determines in a certain sense the rate and size at which the jumps of the process occur (for a Brownian motion, there are no jumps, and $\Pi\equiv 0$.) $\Pi$ is in turn determined by its positive and negative tails, which for the stable laws are of the form
\be\label{tailst}
\pibar^+(x)= \frac{c_+}{x^\alpha}, \ 
\pibar^-(x)= \frac{c_-}{x^\alpha}, \ x>0,
\ee
where $c_\pm\ge 0$ are constants with $c=c_++c_->0$,
and $\alpha\in(0,2)$ is the index of the stable law.
We refer to
\be\label{btail}
\pibar(x)= \pibar^+(x)+\pibar^-(x)= \frac{c}{x^\alpha},\ x>0, 
\ee
simply as the {\it tail} of the measure. An important point is that $\pibar$ possesses a singularity at $x=0$ of order $\alpha$; this corresponds to a concentration of infinitely many small jumps w.p.1\footnote{``w.p.1"= ``with probability 1", or almost surely}
 in any finite interval. One can see the possible relevance of this property for modelling in the thin wafer situation.

  For further properties of stable processes we refer to \cite{satoa} and \cite{st}.

\section{Passage Times out of Rectangular or Curvilinear Regions}\label{s5} 
In this section we consider the passage times $T_\kappa(r)$  of a general L\'evy process defined in \eqref{defT}, but including special reference to the stable processes.
To begin with we note a distinct difference between the cases  $\kappa=0$ and  $\kappa>0$.
When   $\kappa=0$ and $r>0$ the region is a rectangular (slab or wafer) of width $2r$ not including (0,0). 
But when   $\kappa>0$ the curvilinear region shrinks to 0 at small times.
In the second case, there is the possibility that the continuous time process might jump out of the region instantaneously, even for a fixed value of $r>0$. At another extreme, there is the possibility that the process never gets out of the region at all; in this case we formally set $T_\kappa(r)=\infty$. The latter is never the case when $\kappa=0$ but can occur when $\kappa>0$ since the width of the intervals then expands indefinitely.
We need to keep these cases in mind when analysing the passage times. In general they hinge on the relative sizes of the two parameters $\kappa$ and $\alpha$, as we now discuss.

\bigskip\noindent{\it Instantaneous Exit}\

To rule out the possibility of an  instantaneous exit
we have the following equivalences, true for any L\'evy process $X$, which can be obtained by elementary calculations. For any $r>0$, 
\be\label{A}
P(T_\kappa(r)>0)>0
\iff
P(T_\kappa(r)>0)=1
\iff 
\liminf_{t\dto 0}\frac{|X_t|}{t^\kappa}\le r,\
{\rm w.p.1}.
\ee
Instantaneous exit is thus ruled out when the righthand relation in \eqref{A} holds. 
We can get information on this from \cite{bdm}, Theorem 2.1. From it we deduce that
for $\kappa<1/2$,   $\lim_{t\dto 0}|X_t|/t^\kappa=0$ w.p.1 for any L\'evy process, so the righthand relation in \eqref{A} always holds in this case;  while for  $\kappa>1/2$, 
 $\lim_{t\dto 0}|X_t|/t^\kappa=0$ w.p.1 if and only if
 (iff) $\int_0^1\pibar(x^\kappa)\rmd x$ is finite.
 Referring to \eqref{btail}, we can conclude:
{\it Instantaneous exit is ruled out for the stable process $S^{(\alpha)}$ if $0\le \kappa<1/\alpha$.}
 The interpretation is that as long as the singularity in $\Pi$ at 0 is not too steep ($\alpha<1/\kappa$),
 $S_t^{(\alpha)}$ takes a nontrivial time to escape from $R_\kappa$.
 
  \bigskip\noindent{\it  Exit in Finite Time}\

When $T_\kappa(r)=0$, exit obviously occurs in finite time,  but when  $T_\kappa(r)>0$ there is the possibility that $T_\kappa(r)=\infty$, which is to say that w.p.1 $X$ never escapes from the region in finite time. This situation can be analysed as follows.
It is convenient to define $\Xbar^\kappa_t:=\sup_{0<s\le t}|X_s|/s^\kappa$, $t>0$, with
 $\Xbar^\kappa_\infty:=\sup_{t>0}|X_t|/t^\kappa$.
Then for $t>0$
\be\label{subs}
\{T_\kappa(r)>t\} \subseteq
\{\Xbar^\kappa_t\le r\} \subseteq\{T_\kappa(r)\ge t\}.
\ee
In the Appendix Section \ref{s1.4} some equivalences are proved for $P(T_\kappa(r)<\infty)=1$, and it is shown that
this holds iff $P(|X_t|>rt^\kappa\ {\rm for\ some}\ t>0)=1$.  Sufficient for this is that
\be\label{D}
\limsup_{t\dto 0}\frac{|X_t|}{t^\kappa}>r,\
{\rm or}\ \limsup_{t\to\infty}\frac{|X_t|}{t^\kappa}>r,\ {\rm w.p.1}.
\ee
\eqref{D} describes escape of $X_t$ from $R_k$ infinitely often w.p.1 either at small times or at large times. But it is easier to deal with than
the equivalences in \eqref{eqch} of the Appendix.

Theorem 2.1 of \cite{bdm} can be used again for information on the lefthand condition in \eqref{D}.
 For $\kappa<1/2$,   it never holds, 
 while for  $\kappa>1/2$, we can deduce that
 $\limsup_{t\dto 0}|X_t|/t^\kappa>0$ w.p.1
 iff $\limsup_{t\dto 0}|X_t|/t^\kappa=\infty$ w.p.1 iff
  $\int_0^1\pibar(x^\kappa)\rmd x$ is infinite.
 For the stable process $S_t^{(\alpha)}$, the latter is the case iff $0\le \kappa<1/\alpha$.

For the righthand condition in \eqref{D}, we refer to
\cite{dm3}, which gives 
 $\limsup_{t\to\infty}|X_t|/t^\kappa=\infty$ w.p.1
except when 
$\kappa\ge 1$ or $1/2<\kappa<1, E|X_1|<\infty$ and $EX_1=0$.
In these latter cases, 
 $\limsup_{t\to\infty}|X_t|/t^\kappa=\infty$ w.p.1
 iff 
  $\int_1^\infty\pibar(x^\kappa)\rmd x$ is infinite.
    Specialising to the stable process $S_t^{(\alpha)}$,   the integral diverges iff $\kappa<1/\alpha$.
  Putting together both sides of \eqref{D}, we conclude:  $S_t^{(\alpha)}$ 
  {\it  escapes from $R_\kappa$ in finite time w.p.1 for any $\kappa\ge 0$ and $0<\alpha<2$.  }

\bigskip\noindent{\it Moments of the Passage Times}\

Restricting ourselves at first to the discrete time situation, 
we can deduce from Prop. 2.2 of Doney and Griffin (2003)  that $ET_\kappa^S(r)$ is finite  for any $r>0$ when $0\le \kappa<1/2$. 
The same property carries over to the continuous time version.

\section{Asymptotics}\label{s6}
The results  in the previous section apply for each fixed $r>0$.  In a physical situation we expect the process will only approximately behave as a stable process. The approximation may be close (and desired to be close) at small times (corresponding to small slab widths $r$ for the passage times) or  it may be close at large times (corresponding to large widths $r$).
If a L\'evy flight (discrete time heavy tailed random walk) is assumed for the motion, \eqref{doa} and \eqref{dofna} describe approximation by a continuous time stable process. 
If the underlying model is taken to be a general (continuous time) L\'evy process, it may also be well described by a stable process for small or large time.  

In either situation, because of the functional laws, the passage times in continuous time (defined by \eqref{defT})  or in discrete time (defined by \eqref{defTS}), 
converge after appropriate rescaling to the corresponding passage times of the approximating stable process either as $r\dto 0$ or as $r\to\infty$, corresponding to small or large time approximations. These have the properties listed in the previous section, and others
related to the distributions of the passage times.

\subsection{Relative Stability of Passage Times}\label{s8}
``Relative Stability" refers to the convergence in distribution of $T_\kappa(r)$, either as $r\dto 0$ or as $r\to\infty$,
to a nonzero degenerate random variable. Thus, convergence
is in fact in probability to a finite nonzero constant. 
This is a special case of the  convergence in distribution of 
$T_\kappa(r)$. 
Specifically, we have relative stability when $T_\kappa(r)/C(r)\topr 1$,
for some nonstochastic function $C(r)>0$,  where the convergence is in probability as $r\dto 0$ or as $r\to\infty$. 

This kind of behaviour occurs iff the underlying process $X$ is ``relatively stable"  in 
the  sense that $X_t$ converges in probability to a finite nonzero constant,  either as $t\dto 0$ or as $t\to\infty$. 
This property has been thoroughly analysed and conditions for it to occur are listed in \cite{gm2011} and \cite{gm2013}. They hold in particular when  $X$ has bounded variation with positive drift (for $r\dto 0$), or when $X$ has finite, positive mean (for $r\to\infty$). In these cases the norming function $C(r)$ described in the previous paragraph can be taken as  a constant multiple of $r^{1/(1-\kappa)}$, $0\le \kappa<1$. When $\kappa=0$, corresponding to exit from a slab of width $2r$, this means that  $T_\kappa(r)$ behaves approximately linearly  as $r\dto 0$ or as $r\to\infty$.

\subsection{Passage Times Distributions}\label{s7}
Returning to the random walk model defined in \eqref{defS}, 
a relevant large-time functional law for the passage time is in
Doney and Maller (2005). 
It tells us that when \eqref{dofna} holds 
\be\label{Tcon}
\frac{T_\kappa^S(r)}{r^\nu}\todr 
\frac{1}{Y^{1/\nu}}, \ {\rm as}\ r\to \infty,
\ee
where $\nu=1/\alpha-\kappa$ and $\kappa<1/\alpha$. In \eqref{Tcon}, $Y$ has the distribution of 
$\sup_{0<s\le t}|S_t^{(\alpha)}|/t^\kappa$.
The restriction $\kappa<1/\alpha$ ensures that $Y$ is finite w.p.1.
Similar relations  hold for the continuous time $T_\kappa(r)$ defined in \eqref{defT}.

\section{Overview}\label{s9}
The following overview may help to distinguish the various processes discussed in the paper. 
We may have as model:

\bigskip
$\bullet$ a discrete time random walk $(S_n)_{n=1,2,\ldots}$ in the domain of attraction of a non-normal stable law as $n\to\infty$ 
(i.e., the  random walk is a ``L\'evy flight")

$\bullet$ which converges as $n\to\infty$,  after rescaling, to a non-normal stable process  $(S_t^{(\alpha)})_{t>0}$ (via the functional law in \eqref{dofna}), so that
 
$\bullet$ the discrete time passage time in \eqref{defTS} converges to the stable process passage time in \eqref{Talpha}
(again, via the functional law).

\bigskip
Alternatively,  we may have as model:

$\bullet$ a continuous time L\'evy process $(X_t)_{t\ge 0}$ in the domain of attraction of a stable law either as $t\dto 0$ or as $t\to\infty$

$\bullet$ which converges either as $t\dto 0$ or as $t\to\infty$ to a stable process  $(S_t^{(\alpha)})$, $t>0$, (via a corresponding functional law), so that
 
$\bullet$ the continuous time passage time in \eqref{defT} converges to the stable process passage time in \eqref{Talpha}
(via the functional law).

\noindent
 Alternatively,  again, the continuous time model may be the stable process  $(S_t^{(\alpha)})$, itself. 

\bigskip
Having constructed passage times $T_\kappa^S(r)$, $T_\kappa(r)$ or $T_\kappa^{(\alpha)}(r)$ from one of these processes, we can study the limiting behaviour as $r\to\infty$ (for any of $T_\kappa^S(r)$, $T_\kappa(r)$ or $T_\kappa^{(\alpha)}(r)$) or  as $r\dto 0$ (for $T_\kappa(r)$ and $T_\kappa^{(\alpha)}(r)$).

\section{Summary}\label{s10}
There is a vast literature on domains of attraction both for continuous and discrete time processes, and for stable distributions and processes themselves. We have discussed only the 2-sided exit problem as this seems to be relatively newly considered in the applications literature. One-sided exits, i.e., from a semi-infinite domain, have been more commonly studied and applied, for example in the fluctuation theory of random walks and L\'evy processes.  Still, many gaps remain to be filled.

\section{Appendix: Notes and Further References}\label{s1.4}
For exit in finite time we have the following equivalences. Recall the definition
 $\Xbar^\kappa_t:=\sup_{0<s\le t}|X_s|/s^\kappa$, $t>0$, with
 $\Xbar^\kappa_\infty:=\sup_{t>0}|X_t|/t^\kappa$.
Then 
\bea\label{eqch}
&&
P(T_\kappa(r)<\infty)=1
\iff
\lim_{t\to\infty} P(T_\kappa(r)>t)=0
\iff 
\lim_{t\to\infty} P(\Xbar_t^\kappa\le r)=0\iff
\cr &&
P(\Xbar_\infty^\kappa\le r)=0
\iff
P(\sup_{t>0}|X_t|/t^\kappa\le r)=0
\iff
P(|X_t|>rt^\kappa\ {\rm for\ some}\ t>0)=1.\cr &&
\eea
Some information on the rate of divergence of  $ET_\kappa^S(r)$ as $r\to\infty$ is in \cite{DG} and \cite{DG2}.

\bigskip\noindent {\bf Acknowledgements}\ 
This research was partially supported by ARC Grant  DP1092502.

\end{document}